\documentclass[11pt,leqno,fleqn]{article}
\usepackage{amsmath,amssymb}
\newcommand{\dps}{\displaystyle}
\newcommand{\End}{\text{\rm End}}
\newcommand{\Hom}{\text{\rm Hom}}
\newcommand{\sgn}{\text{\rm sgn}}
\newcommand{\GL}{\text{\rm GL}}
\newcommand{\SL}{\text{\rm SL}}
\newcommand{\SO}{\text{\rm SO}}
\newcommand{\Sp}{\text{\rm Sp}}
\newcommand{\T}{^{\sf T}}
\newcommand{\dy}[2]{%
\refstepcounter{equation}%
\label{#1}%
\begin{list}{}{
\topsep 5mm
\leftmargin 18mm
\rightmargin 0cm
\itemsep 0mm
\listparindent 0mm
\parsep 0mm
\itemsep 0mm
\labelsep 0mm
\labelwidth 18mm
}%
\item[\rm (\theequation)\hfill]
#2
\end{list}%
}
\newcommand{\dyz}[1]{%
\refstepcounter{equation}%
\begin{list}{}{
\topsep 5mm
\leftmargin 18mm
\rightmargin 0cm
\itemsep 0mm
\listparindent 0mm
\parsep 0mm
\itemsep 0mm
\labelsep 0mm
\labelwidth 18mm
}%
\item[\rm (\theequation)\hfill]
#1
\end{list}%
}
\newcommand{\dyyz}[1]{\dyz{\raggedright$\dps#1$}}
\newcommand{\dyy}[2]{\dy{#1}{\raggedright$\dps#2$}}
\newcommand{\de}[2]{\dy{#1}{\raggedright$\displaystyle #2 $}}
\newcommand{\dez}[1]{\dyz{\raggedright$\displaystyle #1 $}}
\newcounter{stelling}
\newcommand{\thm}[2]{\setcounter{gevolg}{0}\setcounter{claim}{0}\refstepcounter{stelling}\vspace{4mm}\noindent{\bf Theorem \thestelling.}\label{#1}{\it #2}}
\newcounter{gevolg}[stelling]
\newcommand{\cor}[2]{\refstepcounter{gevolg}\setcounter{claim}{0}\vspace{4mm}\noindent{\bf Corollary \thestelling.\thegevolg.}\label{#1}{\it #2}}
\newcommand{\corz}[1]{\refstepcounter{gevolg}\setcounter{claim}{0}\vspace{4mm}\noindent{\bf Corollary \thestelling.\thegevolg.}{\it #1}}
\newcounter{claim}
\newcommand{\cl}[2]{\refstepcounter{claim}\vspace{4mm}\noindent{\em Claim \theclaim.} \label{#1} {\it #2}}
\newcommand{\clz}[1]{\refstepcounter{claim}\vspace{4mm}\noindent{\em Claim \theclaim.}  {\it #1}}
\newcounter{sectie}
\newcommand{\sect}[2]{\refstepcounter{sectie}
\section*{\boldmath \thesectie. #2}%
\label{#1}}
\newcommand{\sectz}[1]{\refstepcounter{sectie}
\section*{\boldmath \thesectie. #1}%
}
\newcommand{\pf}{\vspace{3mm}\noindent{\bf Proof.}\ }
\newcommand{\pfcl}{\vspace{3mm}\noindent{\em Proof.}\ }
\newcommand{\bx}{\hspace*{\fill} \hbox{\hskip 1pt \vrule width 4pt height 8pt depth 1.5pt \hskip 1pt}

\addvspace{4mm}}
\newcommand{\obx}{\hspace*{\fill} \hbox{$\Box$}

\addvspace{4mm}}
\newcommand{\bxx}{\hspace*{\fill} \hbox{\hskip 1pt \vrule width 4pt height 8pt depth 1.5pt \hskip 1pt}}
\newcommand{\rf}[1]{{\rm (\ref{#1})}}
\newcommand{\tr}{\text{\rm tr}}
\newcommand{\AAA}{{\cal A}}
\newcommand{\CC}{{\cal C}}
\newcommand{\EE}{{\cal E}}
\newcommand{\FF}{{\cal F}}
\newcommand{\II}{{\cal I}}
\newcommand{\MM}{{\cal M}}
\renewcommand{\SS}{{\cal S}}
\newcommand{\UU}{{\cal U}}
\newcommand{\linhull}{\text{\rm lin.hull}}
\newcommand{\oC}{{\mathbb{C}}}
\begin{document}

\begin{center}
{\LARGE\bf Tensor subalgebras and First Fundamental Theorems in
invariant theory

}

\end{center}

\begin{center}
{\large
Alexander Schrijver\footnote{ CWI and University of Amsterdam.
Mailing address: CWI, Kruislaan 413, 1098 SJ Amsterdam,
The Netherlands.
Email: lex@cwi.nl.}}

\end{center}

\noindent
{\small{\bf Abstract.}
Let $V=\oC^n$ and let $T:=T(V)\otimes T(V^*)$ be the mixed tensor algebra
over $V$.
We characterize those subsets $A$ of $T$ for which there is a subgroup $G$
of the unitary group $\UU(n)$ such that $A=T^G$.
They are precisely the nondegenerate contraction-closed graded $*$-subalgebras
of $T$.
While the proof makes use of the First Fundamental Theorem for
$\GL(n,\oC)$ (in the sense of Weyl), the characterization has as direct
consequences First
Fundamental Theorems for several subgroups of $\GL(n,\oC)$.
Moreover, a Galois connection between linear algebraic $*$-subgroups of
$\GL(n,\oC)$ and nondegenerate contraction-closed $*$-subalgebras of
$T$ is derived.

}

\sectz{Introduction}

Let $V=\oC^n$.
Denote, as usual,
\dyz{
$\dps T(V):=\bigoplus_{k=0}^{\infty}V^{\otimes k}$ 
and
$\dps T(V^*):=\bigoplus_{k=0}^{\infty}{V^*}^{\otimes k}$,
}
where $V^{\otimes k}$ and $V^{*\otimes k}$ denote the tensor product
of $k$ copies of $V$ and $V^*$ respectively.
Set
\dyyz{
T:=T(V)\otimes T(V^*)\cong
\bigoplus_{k,l=0}^{\infty}V^{\otimes k}\otimes {V^*}^{\otimes l}.
}
This is the {\em mixed tensor algebra} over $V$ (cf.\ [3]).
(The multiplication is the usual tensor product of the rings $T(V)$ and
$T(V^*)$, governed by the rule
$(x\otimes y)\otimes(x'\otimes y')=(x\otimes x')\otimes(y\otimes y')$ for
$x,x'\in T(V)$ and $y,y'\in T(V^*)$.)

For any $U\in\GL(n,\oC)$, let $z\mapsto z^U$ be the action
of $U$ on $T$, which is the unique algebra
endomorphism satisfying $x^U=Ux$ and $y^U=(U\T)^{-1}y$ for $x\in V$ and
$y\in V^*$.
For any $G\subseteq\GL(n,\oC)$ and $A\subseteq T$, denote
\dyz{
$A^G:=\{z\in A\mid z^U=z$ for all $U\in G\}$ and
$G^A:=\{U\in G\mid z^U=z$ for all $z\in A\}$.
}

In this paper we characterize those subsets $A$ of $T$
for which there exists a subgroup $G$ of the unitary group $\UU(n)$
such that $A=T^G$.
They turn out to be precisely the graded $*$-subalgebras of
$T$ that are nondegenerate and
contraction-closed (for definitions, see Section \ref{17ja06a}).
Our proof is based on the Stone-Weierstrass theorem, the
First Fundamental Theorem (in the sense of Weyl [4]) for
$\GL(n,\oC)$, and the existence of a Haar measure on $\UU(n)$.

As consequences we derive the First Fundamental Theorem
for a number of subgroups of $\GL(n,\oC)$.
Indeed, our theorem directly implies that if some subgroup $G$ of
$\GL(n,\oC)$ satisfies $G=G^S$ for some subset $S$ of
$T$ with $S=S^*$, then $T^G$ is equal to the
smallest nondegenerate contraction-closed graded subalgebra
of $T$ containing $S$.
This often directly yields a spanning set for
$(V^{\otimes k}\otimes V^{*\otimes l})^G$ for all $k,l$.
That is, it implies a First Fundamental Theorem for $G$
(the tensor version, which is equivalent to the polynomial version
--- cf.\ Goodman and Wallach [2] Section 4.2.3).
We describe this in Section \ref{15ja06b}.

A subgroup $G$ of $\GL(n,\oC)$ is called a {\em $*$-subgroup} if
$G=G^*:=\{U^*\mid U\in G\}$.
In Section \ref{15ja06c} we show that for any $*$-subgroup $G$ of
$\GL(n,\oC)$, if we set $A:=T^G$, then
$\GL(n,\oC)^A$ is equal to the smallest linear algebraic subgroup of
$\GL(n,\oC)$ containing $G$.
Together with the characterization above, this implies a Galois connection
between linear algebraic $*$-subgroups of $\GL(n,\oC)$ and nondegenerate
contraction-closed graded $*$-subalgebras of $T$.

\sect{17ja06a}{Preliminaries}

For any $A\subseteq T$ and $k,l\geq 0$, denote
\dez{
A_{k,l}:=A\cap (V^{\otimes k}\otimes V^{*\otimes l}).
}
A subalgebra $A$ of $T$ is called {\em graded} if
$A=\bigoplus_{k,l=0}^{\infty}A_{k,l}$.

Fix a basis $\{e_1,\ldots,e_n\}$ of $V$, and the
corresponding dual basis $e_1^*,\ldots,e^*_n$ of $V^*$.
This extends to a unique function $x\mapsto x^*$
on $T$ satisfying
$(x^*)^*=x$, $(\lambda x)^*=\overline\lambda x^*$, $(x+y)^*=x^*+y^*$,
and $(x\otimes y)^*=y^*\otimes x^*$ for all $x,y\in T$
and $\lambda\in\oC$.
A subalgebra $A$ of $T$
is called a {\em $*$-subalgebra} if $A^*=A$.
Note that if $U$ is unitary, then $(z^*)^U=(z^U)^*$ for all
$z\in T$.

Call a subalgebra $A$ of $T$ {\em nondegenerate} if
there is no proper subspace $W$ of $V$ such that
$A\subseteq T(W)\otimes T(W^*)$.
A subalgebra is {\em contraction-closed} if it is closed under
the contraction operators.
We call a tensor $z\in V^{\otimes k}\otimes V^{*\otimes l}$ a
{\em mutation} of a tensor $y\in V^{\otimes k}\otimes V^{*\otimes l}$
if $z$ arises from $y$ by permuting contravariant factors and permuting
covariant factors.

Finally, we denote
\dez{
I:=\sum_{i=1}^n e_i\otimes e_i^*,
}
the {\em identity matrix}.

Most background on tensors and invariant theory can be found
in the book of Goodman and Wallach [2].

\sectz{The characterization}

\thm{12ja06c}{
Let $n\geq 1$ and $A\subseteq T$.
Then there is a subgroup $G$ of $\UU(n)$ with $A=T^G$
if and only if $A$ is a nondegenerate contraction-closed graded
$*$-subalgebra of $T$.
}

\pf
Necessity being direct, we show sufficiency.
Let $A$ be a nondegenerate contraction-closed graded $*$-subalgebra of
$T$.
We first show:

\cl{12ja06a}{
$I\in A$.
}

\pfcl
Consider the elements of $V\otimes V^*$ as elements of $\End(V)$,
or as the corresponding matrices.
Then $A_{1,1}$ is a subalgebra of $\End(V)$, since if $y,z\in A_{1,1}$
then the matrix product $yz$ belongs to $A_{1,1}$, as it is a contraction
of $y\otimes z$.
For each $y\in A_{1,1}$, let $C_y$ be the column space of $y$.
That is,
\dez{
C_y:=\{yv\mid v\in V\},
}
where $yv$ denotes the product of matrix $y$ and vector $v$.
Elementary matrix theory tells us that for all matrices $y,z$ we have
$C_y=C_{yy^*}$ and $C_y+C_z=C_{yy^*+zz^*}$.
Hence the union of the $C_y$ over all
$y\in A_{1,1}$ is a subspace $W$ of $V$, and $W=C_z$ for some
$z\in A_{1,1}$.
Then
\de{15me05b}{
A\subseteq T(W)\otimes T(W^*).
}
To prove this, we can assume that $W\cap\{e_1,\ldots,e_n\}$ is
a basis of $W$, say it is $\{e_1,\ldots,e_m\}$.
Express any $x\in A_{k,l}$ in the basis of
$V^{\otimes k}\otimes V^{*\otimes l}$ consisting of all tensors
\de{14ja06a}{
e_{i_1}\otimes\cdots\otimes e_{i_k}\otimes e^*_{j_1}\otimes\cdots\otimes e^*_{j_l}
}
with $i_1,\ldots,i_k,j_1,\ldots,j_l\in\{1,\ldots,n\}$.
If $x\not\in W^{\otimes k}\otimes W^{*\otimes l}$, then we may assume
(by the fact that $A=A^*$) that
$x$ uses a basis element \rf{14ja06a} with $i_t>m$ for some $t\in\{1,\ldots,k\}$.
Then there is a contraction of $x\otimes x^*$ to an element $y$ of
$A_{1,1}$ which uses $e_{i_t}\otimes e^*_{i_t}$.
Hence $C_y$ contains an element using $e_{i_t}$, contradicting
the fact that $C_y\subseteq W$.
This proves \rf{15me05b}.

As $A$ is nondegenerate, \rf{15me05b} implies $W=V$.
So $z$ has rank $n$.
Hence $I$ is a linear combination of the matrices
$zz^*,(zz^*)^2,(zz^*)^3,\ldots$.
Therefore, $I\in A_{1,1}$.
\obx

Since $A$ is contraction-closed, Claim \ref{12ja06a} implies that
if $x\in A$, then also each mutation of $x$ belongs to $A$ (since each mutation
of $x$ can be obtained by a series of contractions from $x\otimes I^{\otimes k}$
for some $k$).
Moreover, by the First Fundamental Theorem for $\GL(n,\oC)$, Claim \ref{12ja06a}
implies
\de{15ja06f}{
T^{\UU(n)}\subseteq A.
}

Define $G:=\UU(n)^A$.
To prove the theorem, it suffices to show
\dez{
A=T^G.
}
Here $\subseteq$ is direct.
To show the reverse inclusion, choose
$a\in (V^{\otimes k_0}\otimes V^{*\otimes l_0})^G$ for some $k_0,l_0$.
We are done if we have shown that $a\in A$.

Let $X:=\UU(n)/G$, the set of left cosets of $G$.
Let $\langle x,y\rangle$ be the inner product on $T$
determined by the basis $e_1,\ldots,e_n$.
(So $\langle x,y\rangle=0$ if $x\in V^{\otimes k}\otimes V^{*\otimes l}$
and $y\in V^{\otimes k'}\otimes V^{*\otimes l'}$ with
$(k,l)\neq (k',l')$.)

For $y\in T$ and $z\in A$, define a function
$\phi_{y,z}:X\to\oC$ by
\dez{
\phi_{y,z}(UG):=\langle y,z^U\rangle
}
for $U\in \UU(n)$.
This is well-defined, since if $UG=U'G$, then $U^{-1}U'\in G$, hence
$z^{U^{-1}U'}=z$, and therefore $z^U=z^{U'}$.

Let
\dez{
\FF:=\linhull\{\phi_{y,z}\mid y\in T, z\in A\}.
}
We will show that $\overline{\FF}=\CC(X)$ (with respect to the
sup-norm on $\CC(X)$), by applying the Stone-Weierstrass theorem
(cf.\ for instance [1] Corollary 18.10).

First, $\FF$ is a subalgebra of $\CC(X)$ (with respect to
pointwise multiplication).
For let
$y\in V^{\otimes k}\otimes V^{*\otimes l}$,
$z\in A_{k,l}$,
$y'\in V^{\otimes k'}\otimes V^{*\otimes l'}$,
and $z'\in A_{k',l'}$.
Then for each $U\in \UU(n)$:
\dyyz{
\phi_{y,z}(UG)\phi_{y',z'}(UG)=\phi_{y\otimes y',z\otimes z'}(UG).
}
Moreover, $\FF$ is self-conjugate: if $\phi\in\FF$, also
$\overline{\phi}\in\FF$ (as $\overline{\phi_{y,z}}=\phi_{y^*,z^*}$).

Finally,
$\FF$ is strongly separating.
Indeed, for $U,U'\in \UU(n)$ with
$UG\neq U'G$ there exist $y\in T$ and $z\in A$ such that
$\phi_{y,z}(UG)\neq\phi_{y,z}(U'G)$ (as by definition of $G$, there
is a $z\in A$ with $z^{U^{-1}U'}\neq z$ (as $U^{-1}U'\not\in G$),
so $z^U\neq z^{U'}$).
Moreover, for $U\in \UU(n)$ there exist
$y\in T$ and $z\in A$ with $\phi_{y,z}(UG)\neq 0$,
as $A$ contains a nonzero element $z$, hence $z^U\neq 0$, and so
$\phi_{y,z}(UG)\neq 0$ for some $y$.

As $X$ is a compact space, the Stone-Weierstrass theorem
gives $\overline{\FF}=\CC(X)$.
This implies:

\clz{
There exist $n_1,n_2,\ldots$ such that for each $i=1,2,\ldots$
there are $y_{i,j}\in T$, $z_{i,j}\in A$, and
$x_{i,j}\in V^{\otimes k_0}\otimes V^{*\otimes l_0}$
(for $j=1,\ldots,n_i$)
with the property that for each $U\in \UU(n)$:
\de{4fe06a}{
a^U=\lim_{i\to\infty}\sum_{j=1}^{n_i}\phi_{y_{i,j},z_{i,j}}(UG)x_{i,j},
}
where the limit is uniform in $U$.
}

\pfcl
Let $b_1,\ldots,b_t$ be any basis of
$V^{\otimes k_0}\otimes V^{*\otimes l_0}$ and
let $c_1,\ldots,c_t$ be the corresponding dual basis.
(So $t=n^{k_0+l_0}$.)

Consider any fixed $h\in\{1,\ldots,t\}$.
The function $UG\mapsto c_h(a^U)$ is well-defined (as $a^U=a^{U'}$ if
$UG=U'G$, since $a\in T^G$) and belongs to $\CC(X)$.
So the Stone-Weierstrass theorem implies that the function
$UG\mapsto c_h(a^U)$ is the uniform limit of sums of functions
$\phi_{y,z}$.
That is, there exist $y_{i,j,h}\in T$ and $z_{i,j,h}\in A$
(for $i=1,2,\ldots$ and $j=1,\ldots,n_{i,h}$, for some $n_{i,h}$)
such that for each $U\in\UU(n)$:
\dez{
c_h(a^U)=
\lim_{i\to\infty}
\sum_{j=1}^{n_{i,h}}\phi_{y_{i,j,h},z_{i,j,h}}(UG),
}
where the limit is uniform in $U$.

As this holds for each $h=1,\ldots,t$, we have
\dyyz{
a^U=
\sum_{h=1}^tc_h(a^U)b_h
=
\lim_{i\to\infty}
\sum_{h=1}^t
\sum_{j=1}^{n_{i,h}}\phi_{y_{i,j,h},z_{i,j,h}}(UG)b_h.
}
Combining $j,h$ to one new index $j$ gives \rf{4fe06a}.
\obx

\rf{4fe06a} is, by definition of $\phi_{y,z}$, equivalent to:
\dez{
a^U=\lim_{i\to\infty}\sum_{j=1}^{n_i}\langle y_{i,j},z^U_{i,j}\rangle x_{i,j}.
}
Equivalently, for each $U\in \UU(n)$:
\dez{
a=\lim_{i\to\infty}\sum_{j=1}^{n_i}\langle y^U_{i,j},z_{i,j}\rangle x_{i,j}^U.
}
Consider now $x_{i,j}^U\otimes y_{i,j}^{*U}$.
It belongs to
$(V^{\otimes k_0}\otimes V^{*\otimes l_0})\otimes T^*$,
and hence it can be considered as element of
$\Hom(T,V^{\otimes k_0}\otimes V^{*\otimes l_0})$.
So we can view $x_{i,j}^U\otimes y_{i,j}^{*U}$ as matrix.
Then
\dez{
\langle y_{i,j}^U,z_{i,j}\rangle x_{i,j}^U=
(x_{i,j}^U\otimes y_{i,j}^{*U})z_{i,j},
}
the latter being a product of a matrix and a vector.
Now define for all $i,j$:
\dez{
M_{i,j}:=\int_{\UU(n)}(x^U_{i,j}\otimes y^{*U}_{i,j})d\mu(U),
}
where $\mu$ is the normalized invariant Haar measure on $\UU(n)$.
So $M_{i,j}\in (V^{\otimes k_0}\otimes V^{*\otimes l_0})\otimes T^*$.
Since $M_{i,j}$ is $\UU(n)$-invariant,
\rf{15ja06f} implies $M_{i,j}\in A$.
Hence, as $A$ is contraction-closed, $M_{i,j}z_{i,j}\in A$.
Moreover,
\dez{
\lim_{i\to\infty}\sum_{j=1}^{n_i}M_{i,j}z_{i,j}
=
\lim_{i\to\infty}\int_{\UU(n)}\sum_{j=1}^{n_i}(x_{i,j}^U\otimes y_{i,j}^{*U})z_{i,j}d\mu(U)
=
\int_{\UU(n)}\lim_{i\to\infty}\sum_{j=1}^{n_i}(x_{i,j}^U\otimes y_{i,j}^{*U})z_{i,j}d\mu(U)
=
\int_{\UU(n)}ad\mu(U)
=a,
}
since the limit is uniform in $U$ (so that we can exchange $\int$ and
$\lim$).
As each $M_{i,j}z_{i,j}$ belongs to $A_{k_0,l_0}$ and as
$A_{k_0,l_0}$ is (topologically) closed (since it is a linear subspace of
$V^{\otimes k_0}\otimes V^{*\otimes l_0}$), it follows that $a\in A$, as
required.
\bx

\sectz{Corollaries}

We formulate a few consequences of Theorem \ref{12ja06c}.

\cor{12ja06d}{
Let $A$ be a nondegenerate contraction-closed graded
$*$-subalgebra of\/ $T$.
Then
\dez{
T^{\UU(n)^A}=A.
}
}

\pf
Here $\supseteq$ is direct.
We prove $\subseteq$.
By Theorem \ref{12ja06c}, $A=T^G$ for some
subgroup $G$ of $\UU(n)$.
Hence $G\subseteq \UU(n)^A$.
Therefore
\dyyz{
T^{\UU(n)^A}
\subseteq
T^G
=A.
\bxx
}

\cor{12ja06e}{
Let $S\subseteq T$  and let $G$ be a $*$-subgroup of
$\GL(n,\oC)$ with $\UU(n)^S\subseteq G\subseteq\GL(n,\oC)^S$.
Then $T^G$ is equal to the smallest
contraction-closed graded $*$-subalgebra of\/
$T$ containing\/ $S\cup\{I\}$.
}

\pf
Let $A$ be the smallest contraction-closed graded
$*$-sub\-al\-ge\-bra of $T$ containing $S\cup\{I\}$.
So $A$ is obtained from $S\cup S^*\cup\{I\}$ by a series of
linear combinations, tensor products, and contractions.
Hence $\UU(n)^S=\UU(n)^A=:H$.
So, as $G\supseteq H$ and as $T^G$ is a contraction-closed
graded $*$-subalgebra containing $S\cup\{I\}$,
\dyy{13ja06a}{
A\subseteq T^G \subseteq
T^H=A,
}
by Corollary \ref{12ja06c}.\ref{12ja06d}.
Therefore, we have equality throughout in \rf{13ja06a}, which proves the
corollary.
\bx

This corollary implies that each contraction-closed
graded $*$-subalgebra $A$ of $T$ is finitely
generated {\em as a contraction-closed algebra}.
That is, there is a finite subset $S$ of $A$ such that all elements
of $A$ can be obtained from $S$
by linear combinations, tensor products, and contractions.

\corz{
Each contraction-closed graded
$*$-subalgebra $A$ of $T$ is finitely generated as
a contraction-closed algebra.
}

\pf
We may assume that $A$ is nondegenerate.
Let $G:=\UU(n)^A$, and for each $z\in A$, let $G_z:=\UU(n)^{\{z\}}$.
So $G=\bigcap_{z\in A}G_z$.
As each $G_z$ is determined by polynomial equations,
by Hilbert's finite basis theorem we know that $G=\bigcap_{z\in S}G_z$
for some finite subset $S$ of $A$.
So $G=\UU(n)^{S}$.
Hence
\dez{
A=T^G=
T^{\UU(n)^{S}}.
}
We can assume that $S^*=S$ (otherwise add $S^*$ to $S$).
Therefore, by Corollary \ref{12ja06c}.\ref{12ja06e}, $A$ is the smallest 
contraction-closed graded sub\-al\-ge\-bra of $T$
containing $S\cup\{I\}$.
\bx

\sect{15ja06b}{Applications to FFT's}

We now apply Theorem \ref{12ja06c} (more precisely, Corollary
\ref{12ja06c}.\ref{12ja06e}) to derive a First Fundamental Theorem
(FFT) in the sense of Weyl [4] for a number of subgroups of
$\GL(n,\oC)$.

\bigskip
\noindent
{\bf FFT for $\SL(n,\oC)=\{U\in\GL(n,\oC)\mid\det U=1\}$} (the special
linear group).
Define $\det\in V^{*\otimes n}$ by
\dez{
\det:=\sum_{\pi\in \SS_n}\sgn(\pi)e^*_{\pi(1)}\otimes\cdots\otimes
e^*_{\pi(n)}.
}
(We can consider $\det$ as element of $(V^{\otimes n})^*$, and
then $\det(x_1\otimes\cdots\otimes x_n)$ is equal to the usual
determinant of the matrix with columns $x_1,\ldots,x_n$.)

One directly checks that $\GL(n,\oC)^{\{\det\}}=\SL(n,\oC)$.
So by Corollary \ref{12ja06c}.\ref{12ja06e}, $T^{\SL(n,\oC)}$ is equal to the smallest
contraction-closed graded subalgebra of $T$ containing $\det$,
$\det^*$, and $I$.
Now $\det\otimes\det^*$ is a linear combination of mutations of
$I^{\otimes n}$ (as it belongs to $T^{\GL(n,\oC)}$).
Hence $T^{\SL(n,\oC)}$ is spanned by
mutations of tensor products of $\det$, $\det^*$, and $I$.

\bigskip
\noindent
{\bf FFT for $\SL_k(n,\oC)=\{U\in\GL(n,\oC)\mid\det U^k=1\}$.}
It is direct to check that
$\GL(n,\oC)^{\{\det^{\otimes k}\}}$
$=\SL_k(n,\oC)$.
So by Corollary \ref{12ja06c}.\ref{12ja06e}, $T^{\SL_k(n,\oC)}$ is equal to the smallest
contraction-closed graded subalgebra of $T$ containing
$\det^{\otimes k}$, $\det^{*\otimes k}$, and $I$.
Again, as $\det^{\otimes k}\otimes\det^{*\otimes k}$ is a linear
combination of mutations of $I^{\otimes nk}$, it follows that
$T^{\SL_k(n,\oC)}$ is spanned by mutations of tensor products of
$\det^{\otimes k}$, $\det^{*\otimes k}$, and $I$.

\bigskip
\noindent
{\bf FFT for $\SS_n(\oC)=$} set of matrices in $\GL(n,\oC)$ with
precisely one nonzero in each row (hence also in each column).
For each $k$, let
\dez{
j_k:=\sum_{i=1}^n(e_i\otimes e^*_i)^{\otimes k}
}
and define
\dez{
J:=\{j_k\mid k\geq 1\}.
}
Then $\GL(n,\oC)^{\{j_2\}}=\SS_n(\oC)$, as one easily checks.
So $T^{\SS_n(\oC)}$ is the smallest contraction-closed graded
subalgebra of $T$ containing $j_2$.
Now the contractions of tensor powers of $j_2$ are precisely the
mutations of tensor products of elements of $J$.
Hence $T^{\SS_n(\oC)}$ is spanned by mutations
of tensor products of elements of $J$.

\bigskip
\noindent
{\bf FFT for $O(n,\oC)=\{U\in\GL(n,\oC)\mid UU\T=I\}$} (the orthogonal
group).
Define
\dez{
f:=\sum_{i=1}^ne_i\otimes e_i.
}
Then $\GL(n,\oC)^{\{f\}}=O(n,\oC)$.
Also, any contraction of $f\otimes f^*$ is equal to $I$.
Hence $T^{O(n,\oC)}$ is spanned by mutations of
tensor products of $f$, $f^*$, and $I$.

\bigskip
In describing the FFT for subgroups of $O(n,\oC)$, it is convenient to
introduce the concept of a `flip'.
A {\em flip} of an element $z\in V^{\otimes k}\otimes V^{*\otimes l}$
is obtained by applying the $\oC$-linear transformation $e_i\mapsto e_i^*$
($i=1,\ldots,n$)
to some (or none) of the factors of $z$, and the reverse
transformation to some (or none) of the factors of $z$.
(So $z$ is also flip of itself.)

Then $f$ and $f^*$ are flips of $I$.
Hence another way of stating the FFT for $O(n,\oC)$ is that
$T^{O(n,\oC)}$ is spanned by mutations of tensor
products of flips of $I$.

\bigskip
\noindent
{\bf FFT for $\SO(n,\oC)=O(n,\oC)\cap\SL(n,\oC)$} (the special
orthogonal group).
Now $\GL(n,\oC)^{\{f,\det\}}$
$=\SO(n,\oC)$.
Hence $T^{\SO(n,\oC)}$ is spanned by mutations of
tensor products of flips of $\det$ and $I$.

\bigskip
\noindent
{\bf FFT for $\SS_n=$} set of $n\times n$ permutation matrices
(the symmetric group).
For each $k$, let
\dez{
h_k:=\sum_{i=1}^ne_i^{\otimes k}
}
and define
\dez{
H:=\{h_k\mid k\geq 1\}.
}
Then $\GL(n,\oC)^H=\SS_n$.
Hence $T^{\SS_n}$ is spanned by mutations of
tensor products of flips of elements of $H$.

A second (but now finite), and more familiar, set of spanning tensors
can be derived from it.
For each $k$, define
\dez{
g_k:=\sum_{i_1,\ldots,i_k}e_{i_1}\otimes\cdots\otimes e_{i_k},
}
where the sum ranges over all distinct $i_1,\ldots,i_k\in\{1,\ldots,n\}$.
(So $g_k=0$ if $k>n$.)
Then
\dyyz{
g_k=
\sum_{f:\{1,\ldots,k\}\to\{1,\ldots,n\}}
\sum_{\pi\in S_k\atop f\circ\pi=f}
\sgn(\pi)e_{f(1)}\otimes\cdots\otimes e_{f(k)}
=
\sum_{\pi\in S_k}
\sgn(\pi)
\sum_{f:\{1,\ldots,k\}\to\{1,\ldots,n\}\atop f\circ\pi=f}
e_{f(1)}\otimes\cdots\otimes e_{f(k)}.
}
Now, in the last expression, for each fixed $\pi\in S_k$,
the sum over $f$ is a mutation of $h_{i_1}\otimes\cdots\otimes h_{i_t}$,
where $i_1,\ldots,i_t$ are the orbit sizes of $\pi$.
Then $h_k$ itself occurs when $\pi$ has precisely one orbit.
As this holds for each $k$, it follows inductively that each $h_k$ is
spanned by mutations of tensor products of $g_0,\ldots,g_k$.
So $T^{\SS_n}$ is spanned by mutations of
tensor products of flips of $g_0,\ldots,g_n$.

\bigskip
\noindent
{\bf FFT for $\SS^{\pm}_n=O(n,\oC)\cap\SS_n(\oC)$}
(so each nonzero entry of any matrix in $\SS^{\pm}_n$ is $\pm 1$).
Let
\dyz{
$H':=\{h_k\mid k$ even, $k\geq 2\}$.
}
Then $\GL(n,\oC)^{H'}=\SS^{\pm}_n$.
Hence $T^{\SS^{\pm}_n}$ is spanned by mutations of
tensor products of flips of elements of $H'$.

As above, one may show that equivalently
$T^{\SS^{\pm}_n}$ is spanned by mutations of
tensor products of flips of
\dez{
\sum_{i_1,\ldots,i_k}e^{\otimes 2}_{i_1}\otimes\cdots\otimes e^{\otimes 2}_{i_k}
}
(for $k=1,\ldots,n$), where the sum ranges over all distinct
$i_1,\ldots,i_k\in\{1,\ldots,n\}$.

\bigskip
\noindent
{\bf FFT for $\AAA_n=\SS_n\cap\SO(n,\oC)$} (the alternating group).
Let $H$ be as above.
Then $\GL(n,\oC)^{H\cup\{\det\}}=\AAA_n$.
Hence $T^{\AAA_n}$ is spanned by mutations of
tensor products of flips of elements of $H$ and of elements
\de{15ja05d}{
\sum_{\pi\in\SS_n}\sgn(\pi)e_{\pi(1)}^{\otimes k_1}\otimes\cdots\otimes
e_{\pi(n)}^{\otimes k_n},
}
ranging over all $k_1,\ldots,k_n\geq 0$.

\bigskip
\noindent
{\bf FFT for $\AAA^{\pm}_n=\SS^{\pm}_n\cap\SO(n,\oC)$.}
Let $H'$ be as above.
Then $\GL(n,\oC)^{H'\cup\{\det\}}=\AAA^{\pm}_n$.
Hence $T^{\AAA^{\pm}_n}$ is spanned by mutations of
tensor products of flips of elements of $H'$ and of elements \rf{15ja05d},
ranging over all odd $k_1,\ldots,k_n\geq 1$.

\bigskip
\noindent
{\bf FFT for $\Sp(n,\oC)=$} set of matrices $U\in\GL(n,\oC)$ with
$USU\T =S$, where
\dez{
S=
\left(
\begin{array}{cc}
0&I_m\\
-I_m&0
\end{array}
\right),
}
for $m:=\frac12n$ (assuming $n$ to be even) (the symplectic group).
Here $I_m$ denotes the $m\times m$ identity matrix.
Define
\dez{
s:=\sum_{i=1}^m(e_i\otimes e_{m+i}-e_{m+i}\otimes e_i).
}
Then $\GL(n,\oC)^{\{s\}}=\Sp(n,\oC)$ (by definition of $\Sp(n,\oC)$).
Any contraction of $s\otimes s^*$ is equal to $\pm I$.
Hence $T^{\Sp(n,\oC)}$ is spanned by tensor
products of $s$, $s^*$, and $I$.

\sect{15ja06c}{A Galois connection}

For any subgroup $G$ of $\GL(n,\oC)$, let $\overline G$ be the
Zariski closure of $G$.
So defining, for any $G\subseteq\GL(n,\oC)$, the ideal $\II_G$ by
\dyz{
$\II_G:=\{p\in\oC[x_{1,1},x_{1,2},\ldots,x_{n,n}]\mid
p(U)=0$ for each $U\in G\}$,
}
then
\dyz{
$\overline G=\{U\in \GL(n,\oC)\mid p(U)=0$ for each $p\in \II_G\}$.
}

\thm{16ja06d}{
For any $*$-subgroup $G$ of $\GL(n,\oC)$:
\dez{
\GL(n,\oC)^{T^G}=\overline G.
}
}

\pf
Define $A:=T^G$ and $H:=\GL(n,\oC)^A$.
So we must prove $H=\overline G$.

Then $H\supseteq\overline G$ is trivial, as $H$ is
linear algebraic, since $H$ is determined by the polynomial
equations $z^U=z$ for $z\in A$.

To prove the reverse inclusion, choose $U_0\not\in\overline G$.
We prove $U_0\not\in H$.
As $U_0\not\in\overline G$, there is a polynomial $p$ such that
$p(U)=0$ for each $U\in G$ while $p(U_0)\neq 0$.
Let $p$ have total degree $m$.
Define for each $U\in \GL(n,\oC)$:
\dez{
\pi(U):= \bigoplus_{k=0}^mU^{\otimes k}.
}
So $\pi$ is a representation of $\GL(n,\oC)$ in the matrix algebra
\dez{
\MM:=
\left(
\bigoplus_{k=0}^mV^{\otimes k}
\right)\otimes\left(
\bigoplus_{l=0}^mV^{*\otimes l}
\right)
\cong
\bigoplus_{k,l=0}^m
\left(
V^{\otimes k}\otimes V^{*\otimes l}
\right)
.
}
Let $\EE$ be the enveloping algebra of $\pi|G$, i.e.,
\dez{
\EE:=\linhull\{\pi(U)\mid U\in G\}.
}
Then from $p$ we can obtain
a matrix $F$ in $\MM$ with $\tr(F\pi(U))=0$ for each
$U\in G$ and $\tr(F\pi(U_0))\neq 0$.
So $\pi(U_0)\not\in\EE$.
Hence, as $\EE$ is a C$*$-algebra, by the double commutant theorem
there exists a matrix $M\in\MM$ commuting with all matrices in
$\EE$ but not with $\pi(U_0)$.
As $M$ commutes with $\pi(U)$ for each $U\in G$, we know
$M^U=M$ for each $U\in G$, and so $M\in A$.
On the other hand, $M$ does not commute with $\pi(U_0)$, and
hence $M^{U_0}\neq M$.
Concluding, $U_0\not\in H$.
\bx

Theorems \ref{12ja06c} and \ref{16ja06d} imply that the relation
$G\leftrightarrow T^G$ gives a
one-to-one relation between the lattice of
linear algebraic $*$-subgroups
$G$ of $\GL(n,\oC)$ and the lattice of nondegenerate
contraction-closed graded $*$-subalgebras of
$T$.
It is a {\em Galois connection}: it reverses inclusion.

This also implies for any linear algebraic $*$-subgroup $G$
of $\GL(n,\oC)$:
\dez{
T^{G\cap\UU(n)}=T^G,
}
since (using Theorems \ref{12ja06c} and \ref{16ja06d}, and denoting
$A:=T^G$)
\dyyz{
T^G=
A=
T^{\UU(n)^A}=
T^{G\cap\UU(n)}.
}

This gives for any linear algebraic $*$-subgroup $G$ of $\GL(n,\oC)$:
\dez{
\overline{G\cap\UU(n)}=G,
}
since
\dyyz{
G=\GL(n,\oC)^{T^G}=
\GL(n,\oC)^{T^{G\cap\UU(n)}}=
\overline{G\cap\UU(n)}.
}
So any linear algebraic $*$-subgroup $G$ of $\GL(n,\oC)$ is determined
by its intersection $G\cap\UU(n)$ with the unitary group.

\bigskip
\noindent
{\em Acknowledgement.}
I thank Laci Lov\'asz for stimulating discussions leading to
conjecturing Theorem \ref{12ja06c}.

\section*{References}\label{REF}
{\small
\begin{itemize}{}{
\setlength{\labelwidth}{4mm}
\setlength{\parsep}{0mm}
\setlength{\itemsep}{1mm}
\setlength{\leftmargin}{5mm}
\setlength{\labelsep}{1mm}
}
\item[\mbox{\rm[1]}] A. Brown, C. Pearcy, 
{\em Introduction to Operator Theory I --- Elements of
Functional Analysis},
Springer, New York, 1977.

\item[\mbox{\rm[2]}] R. Goodman, N.R. Wallach, 
{\em Representations and Invariants of the Classical Groups},
Cambridge University Press, Cambridge, 1998.

\item[\mbox{\rm[3]}] W. Greub, 
{\em Multilinear Algebra --- 2nd Edition},
Springer, New York, 1978.

\item[\mbox{\rm[4]}] H. Weyl, 
{\em The Classical Groups --- Their Invariants and Representations},
Princeton University Press, Princeton, New Jersey, 1946.

\end{itemize}
}

\end{document}